\documentclass[11pt,reqno]{amsart}
\usepackage{amscd,amssymb,amsmath,amsthm}
\usepackage[arrow,matrix]{xy}
\usepackage{graphicx}
\usepackage{cite}
\usepackage{geometry}
\newtheorem{lemma}[subsection]{Lemma}

\numberwithin{equation}{section} 

\newcommand{\bea}{\begin{eqnarray}}
\newcommand{\eea}{\end{eqnarray}}

\begin{document}

\noindent UDC {517.98}

\title{POSITIVE FIXED POINTS OF HAMMERSTEIN'S INTEGRAL OPERATORS WITH DEGENERATE KERNEL}

\author{Yu. Kh. Eshkabilov}

\address{Yu.\ Kh.\ Eshkabilov\\ Karshi State University, 17, Kuchabag st., 180100, Karshi city, Uzbekistan.} \email {yusup62@mail.ru}

\begin{abstract} In this paper we study positive fixed points
of Hammerstein integral operators with degenerate kernel in the
cone of $C[0,1]$. Problem on a number of positive fixed points of
the Hammerstein integral operator leads to the study positive
roots of polynomials with real coefficients. Consider a model on a
Cayley tree with nearest-neighbor interactions and with the set
$[0,1]$ of spin values. The uniqueness translational-invariant
Gibbs measures for the given model is proved.
\end{abstract}
\maketitle


\section{Introduction} \label{sec:intro}

It is well known that integral equations have wide applications in
engineering, mechanics, physics, economics, optimization,
vehicular traffic, biology, queuing theory and so on (see
\cite{gl1994},\cite{ko1972}, \cite{oh1996}, \cite{ab2002},
\cite{kz1984}). The theory of integral equations is rapidly
developing with the help of tools in functional analysis, topology
and fixed point theory. Therefore, many different methods are used
to obtain the solution of the nonlinear integral equation.
Moreover, some methods can be found in Refs. \cite{a2003},
\cite{aee2009}, \cite{aed2005}, \cite{f2005}, \cite{h2008},
\cite{apz2000}, \cite{ak2006}, \cite{b2003}, \cite{m2004} to
discuss and obtain a solution for the Hammerstein integral
equation.The existence of positive solutions of abstract integral
equations of Hammerstein type is discussed in [10].

In this paper we consider an integral operator of Hammerstein type
${{H}_{k}}, k\in \mathbb{N}$ on the cone $C^{+}[0,1]$ of
continuous functions on $[0,1]$:
$$
\left( {{H}_{k}}f \right)\left( t
\right)=\int\limits_{0}^{1}K(t,u)f^{k}du, \eqno{(1.1)}$$

where the kernel $k(t,u)$ is positive continuous function on
$[0,1]^{2}.$

By Theorem 44.8 from [4], the existence of  nontrivial positive
fixed points of the Hammerstein operator (1.1) follows.

In the paper we study a number of positive fixed points of
Hammerstein integral operator with positive degenerate kernel of
the form

$$
\left( {{H}_{k}}f \right)\left( t
\right)=\int\limits_{0}^{1}{\left( {{\varphi }_{1}}\left( t
\right){{\psi }_{1}}\left( u \right)+{{\varphi }_{2}}\left( t
\right){{\psi }_{2}}\left( u \right) \right)}{{f}^{k}}\left( u
\right)du. \eqno{(1.2)}$$

Here ${{\varphi }_{1}}\left( t \right),\,\,\,\,{{\varphi
}_{2}}\left( t \right), \,\  {{\psi }_{1}}\left( t
\right),\,\,\,\,{{\psi }_{2}}\left( t \right)$ are given functions
 from the $C_0^{+}[0,1]=C^{+}[0,1] \setminus \{\theta\}.$

In the case $k=2$  considered the Hammerstein integral operator
$H_k$ with degenerate kernel \cite{enh2016}  and theorems about a
number of positive fixed points of the Hammerstein integral
operator $H_2$ are proved.

The plan of this paper is as follows. In the second section we
study solvability the Hammerstein integral equation with
degenerate kernel in the cone $C_{+}[0,1]:$

$$H_kf=f, k>1. \eqno{(1.3)}$$
In the third section we consider a system of nonlinear algebraic
equations with two unknowns. Problem solvability the system of
nonlinear algebraic equations leads to study positive roots of a
polynomial with the order $k+1$. In the fourth section given
results applied to study Gibbs measures for models on the Cayley
tree $\Gamma^{k}$ of order $k \in \mathbb{N}$.

\section{Integral equation of Hammerstein's type with degenerate
 kernel} \label{sec:intro}

 In this section the existence of positive solution
of the integral equation of Hammerstein type with degenerate
 kernel are discussed.

 Let $C^{+}[0,1]$ is cone of space continuous functions on
 $[0,1].$  We define $C_0^{+}[0,1]=C^{+}[0,1] \setminus \{\theta\}.$
 Let the following functions ${{\varphi }_{1}}\left( t \right),\,\,\,\,{{\varphi }_{2}}\left( t
\right), \,\  {{\psi }_{1}}\left( t \right),\,\,\,\,{{\psi
}_{2}}\left( t \right)$ belongs to the  $C_0^{+}[0,1].$

Define integral operator of the Hammerstein type ${{H}_{k}}, k\in
\mathbb{N}$ on the space $C[0,1]$:
$$
\left( {{H}_{k}}f \right)\left( t
\right)=\int\limits_{0}^{1}{\left( {{\varphi }_{1}}\left( t
\right){{\psi }_{1}}\left( u \right)+{{\varphi }_{2}}\left( t
\right){{\psi }_{2}}\left( u \right) \right)}{{f}^{k}}\left( u
\right)du.$$

We study the integral equation for fixed points of the
Hammerstein's operator ${H}_{k}$:

 $${H}_{k}f=f, \,\,\,\  f \in C_0^{+}[0,1]\eqno{(2.1)}.$$

We define positive numbers ${{a}_{i}}$ and $ {{b}_{i}}$:

${{a}_{i}}=C_{k}^{i}\int\limits_{0}^{1}{{{\psi }_{1}}\left( u
\right)\varphi _{1}^{k-i}\left( u \right)\varphi _{2}^{i}\left( u
\right)du},$ \,\,\,\
 ${{b}_{i}}=C_{k}^{i}\int\limits_{0}^{1}{{{\psi }_{2}}\left( u \right)\varphi _{1}^{k-i}\left( u \right)\varphi _{2}^{i}\left( u \right)du}$,

where $i=\{0,1,2,...,k\}$.

Consider a map $Q_k$ on the two dimensional space
${{\mathbb{R}}^{2}}:$

\[Q_k\left( x,y \right)=\left(
\sum\limits_{i=0}^{k}{{{a}_{i}}{{x}^{k-i}}{{y}^{i}}, \,\
\sum\limits_{i=0}^{k}{{{b}_{i}}{{x}^{k-i}}{{y}^{i}}}} \right).\]

We denote a number of positive fixed points of the operator $T$ by
$N_{+}^{fix}\left(T\right).$

\textbf{Lemma 2.1.} Let $k\geq2$. The Hammerstein's operator
${{H}_{k}}$ has nontrivial positive fixed point iff the map $Q_k$
has nontrivial positive fixed point, moreover $N_{+}^{fix}\left(
{{H}_{k}} \right)=N_{+}^{fix}\left( Q_k \right)$.

\textit{Proof}. a) Let $f\left( t \right)\in C_{0}^{+}\left[ 0,1
\right]$ is nontrivial positive fixed point of the Hammerstein's
operator ${{H}_{k}}$. We introduce the notations
 $${{c}_{1}}=\int\limits_{0}^{1}{{{\psi }_{1}}\left( u \right){{f}^{k}}\left( u \right)}du, \eqno{(2.2)}$$

$${{c}_{2}}=\int\limits_{0}^{1}{{{\psi }_{2}}\left( u \right){{f}^{k}}\left( u \right)}du. \eqno{(2.3)}$$
From the equality $ {{H}_{k}}f=f$ for the fixed point $f$ we have:
$f(t)=c_1\varphi_{1}(t)+c_{2}\varphi_{2}(t).$  Clearly, that
${{c}_{1}}>0,\,\,{{c}_{2}}>0$,  i.e. $( {c}_{1},{c}_{2} ) \in
\mathbb{R}_{2}^{>}$. By the equalities  (2.2) and (2.3) for the
parameters ${{c}_{1}},\,\,{{c}_{2}}$ we obtain the following
equalities
$$ {{c}_{1}}={{a}_{0}}{{c}_{1}}^{k}+{{a}_{1}}{{c}_{1}}^{k-1}{{c}_{2}}+{{a}_{2}}{{c}_{1}}^{k-2}{{c}_{2}}^{2}+\cdots +{{a}_{k-1}}{{c}_{1}}{{c}_{2}}^{k-1}+{{a}_{k}}{{c}_{2}}^{k},
$$
$$
 {{c}_{2}}={{b}_{0}}{{c}_{1}}^{k}+{{b}_{1}}{{c}_{1}}^{k-1}{{c}_{2}}+{{b}_{2}}{{c}_{1}}^{k-2}{{c}_{2}}^{2}+\cdots +{{b}_{k-1}}{{c}_{1}}{{c}_{2}}^{k-1}+{{b}_{k}}{{c}_{2}}^{k}.\\
$$

It means the point
 $( {{c}_{1}},\,{{c}_{2}})$ is fixed point of the map $Q_k.$

b) Let  $\omega=( {{x}_{0}},{{y}_{0}})$ is nontrivial positive
fixed point of the  map $Q_k,$  i.e. $\omega \in
\mathbb{R}_{2}^{+}\backslash \{\theta \}$ and $Q_k\omega=\omega.$
Then
$$ {{a}_{0}}{{x}_{0}}^{k}+{{a}_{1}}{{x}_{0}}^{k-1}{{y}_{0}}+{{a}_{2}}{{x}_{0}}^{k-2}{{y}_{0}}^{2}+\cdots +{{a}_{k-1}}{{x}_{0}}{{y}_{0}}^{k-1}+{{a}_{k}}{{y}_{0}}^{k}={{x}_{0}},
 $$
$$
{{b}_{0}}{{x}_{0}}^{k}+{{b}_{1}}{{x}_{0}}^{k-1}{{y}_{0}}+{{b}_{2}}{{x}_{0}}^{k-2}{{y}_{0}}^{2}+\cdots
+{{b}_{k-1}}{{x}_{0}}{{y}_{0}}^{k-1}+{{b}_{k}}{{y}_{0}}^{k}={{y}_{0}}.
$$
   Using these equalities, we can verify that the function

$${{f}_{0}}(t)={{x}_{0}}{{\varphi }_{1}}(t)+{{y}_{0}}{{\varphi }_{2}}(t)$$
is fixed point of the integral operator ${{H}_{k}}.$\\

\textbf{Theorem 2.1.} Let $k \geq 2.$ The number of postive
solutions of nonlinear integral equation of the Hammerstein type
(2.1) is equal to the number of positive root of the following
polynomial:

$$P_{k+1}(\xi)={{a}_{k}}{{\xi }^{k+1}}+\sum\limits_{i=0}^{k-1}{\left( {{a}_{k-1-i}}-{{b}_{k-i}} \right){{\xi }^{k-i}}}-{{b}_{0}}.
\eqno{(2.4)}$$\\

Proof of the Theorem 2.1 it follows from the Lemma 3.1 and Lemma
3.2 in the section 3.

 By the Theorem 2.1 and by the Descartes rule  for a
number of positive roots of polynomials with real coefficients
(see \cite{p2000} pp. 27-29) it follows \\

 \textbf{Theorem 2.2.} Let $k \geq 2.$

(a) the nonlinear integral equation (2.1) at last one positive
solution;

(b) if for some index  ${{i}_{0}}\in \left\{ 1,\ldots ,k \right\}$
the conditions ${{a}_{i-1}}-{{b}_{i}}\le 0,\,\,1\le i\le
{{i}_{0}}$ and ${{a}_{i-1}}-{{b}_{i}}\ge 0,\,\,{{i}_{0}}<i\le k$
hold then the integral equation (2.1) has unique positive
solution;

(c)  for the number $N_{+}^{fix}\left( {{H}_{k}} \right)$ of the
positive solutions of the nonlinear integral equation (2.1) the
inequality $N_{+}^{fix}\left(
{{H}_{k}} \right)\leq k+1$ holds.\\
Put
$$d_i=a_{i-1}-b_i, \,\,\ i\in \{1, 2, \cdot\cdot\cdot, k\}.$$\\

 \textbf{Corollary 2.1.} Let $k \geq 2.$

a) If for the numbers $d_1, d_2, \cdot\cdot\cdot, d_k$ the
following relations are true:

$$d_1 \leq d_2 \leq \cdot\cdot\cdot \leq d_k$$

then the integral equation (2.1) has unique positive solution;

b) If for the numbers $d_1, d_2, \cdot\cdot\cdot, d_k$ the
following relations are true:

$$d_1\geq d_2 \geq \cdot\cdot\cdot \geq d_k$$

then  for the number $N_{+}^{fix}\left( {{H}_{k}} \right)$ of the
positive solutions of the nonlinear integral equation (2.1) the
inequality $N_{+}^{fix}\left(
{{H}_{k}} \right)\leq 3$ is hold.\\

\textbf{Corollary 2.2.} Let $k \geq 2.$ If there exists two
positive number $\xi_1, \,\ \xi_2$ with $\xi_1<\xi_2$ such that,
$P_{k+1}(\xi_1)>0$ and $P_{k+1}(\xi_2)\leq0$ then
$N_{+}^{fix}\left( {{H}_{k}} \right)\geq
2.$\\

Results of the theorems applied to  study Gibbs measures
\cite{ehr2013}, \cite{ehr2012}, \cite{er2010} for the models  on
the Cayley tree $\Gamma^{2}$.

\section{System of nonlinear algebraic equations
with two unknowns} \label{sec:intro}

In this section we study solvability of system of nonlinear
algebraic equations with two unknowns. We consider the following
system of nonlinear algebraic equations with unknowns $x, y \in
\mathbb{R}:$
$$
\left\lbrace
\begin{array}{ll}
\sum\limits_{i=0}^{k}{a_{i}x^{k-i}y^{i}}=x,\\
\sum\limits_{i=0}^{k}{b_{i}x^{k-i}y^{i}}=y.
\end{array}  \right. \eqno{(3.1)}
$$

\textbf{Lemma 3.1.} If the point $\left( {{x}_{0}},{{y}_{0}}
\right) \in \mathbb{R}_{2}^{>} $ is solution of the system of
nonlinear algebraic equations (3.1) then the number ${{\xi
}_{0}}=\frac{{{y}_{0}}}{{{x}_{0}}}$ is root of the following
polynomial:
$${{a}_{k}}{{\xi }^{k+1}}+\sum\limits_{i=0}^{k-1}{\left( {{a}_{k-1-i}}-{{b}_{k-i}} \right){{\xi }^{k-i}}}-{{b}_{0}}=0.  \eqno{(3.2)}$$

Proof. Let $\left( {{x}_{0}},{{y}_{0}} \right) \in
\mathbb{R}_{2}^{>} $ is solution of the system of nonlinear
algebraic equations (3.1). Then we have

$$
\left\lbrace
\begin{array}{ll}
a_{0}x_0^{k}+a_{1}x_0^{k-1}y_0+\ldots+a_{k-1}x_0y_0^{k-1}+a_{k}y_0^{k}=x_0,\\
b_{0}x_0^{k}+b_{1}x_0^{k-1}y_0+\ldots+b_{k-1}x_0y_0^{k-1}+b_{k}y_0^{k}=y_0.
\end{array}\right.
$$

We introduce the notation ${\xi}=\frac{{y_0}}{{x_0}}$. From here
$y_0=x_0\cdot{\xi}.$ Consequently, we get

$$
\left\lbrace
\begin{array}{ll}
a_{0}x_0^{k}+a_{1}x_0^{k}{\xi}+a_{2}x_0^{k}{\xi}^{2}+\ldots+a_{k-1}x_0^{k}{\xi}^{k-1}+a_{k}x_0^{k}{\xi}^{k}=x_0,\\
b_{0}x_0^{k}+b_{1}x^{k}{\xi}+b_{1}x_0^{k}{\xi}^{2}+\ldots+b_{n-1}x_0^{k}{\xi}^{k-1}+b_{k}x_0^{k}{\xi}^{k}=x_0{\xi}.
\end{array}\right.
$$

Hence,

$$
\left\lbrace
\begin{array}{ll}
x_0^{k}(a_{0}+a_{1}{\xi}+a_{2}{\xi}^{2}+\ldots+a_{k-1}{\xi}^{k-1}+a_{k}{\xi}^{k})=x_0,\\
x_0^{k}(b_{0}+b_{1}{\xi}+b_{1}{\xi}^{2}+\ldots+b_{k-1}{\xi}^{k-1}+b_k{\xi}^{k})=x_0{\xi}.
\end{array}\right.
$$

Consequently

$$\frac{a_{0}+a_{1}{\xi}+a_{2}{\xi}^{2}+\ldots+a_{k-1}{\xi}^{k-1}+a_{k}{\xi}^{k}}{b_{0}+b_{1}{\xi}+b_{1}{\xi}^{2}+\ldots+b_{k-1}{\xi}^{k-1}+b_k{\xi}^{k}}=\frac{1}{\xi}.$$

From the upper equality we obtain

$$a_{k}\xi^{k+1}+(a_{k-1}-b_{k})\xi^{k}+(a_{k-2}-b_{k-1})\xi^{k-1}+\ldots+(a_{0}-b_{1})\xi-b_{0}=0.$$

This completes the proof of the lemma 3.1.

 \textbf{Lemma 3.2.} If the number $\xi_{0}>0$ is root of
the polynomial (3.2) then the point $(x_{0},\xi_{0}x_{0})$ is
solution of the system of nonlinear algebraic equations (3.1),
where

$${{x}_{0}}=\frac{1}{\sqrt[k-1]{\sum\limits_{i=0}^{k}{{{a}_{i}}{{\xi_0 }^{i}}}}}.$$

\textbf{Proof.} Let the positive number $\xi_{0}$ is  root of the
polynomial (3.2) and
\[{{y}_{0}}={{\xi }_{0}}{{x}_{0}},\] i.e.
$${{y}_{0}}=\frac{{{\xi }_{0}}}{\sqrt[k-1]{{{a}_{0}}+{{a}_{1}}{{\xi
}_{0}}+{{a}_{2}}\xi _{0}^{2}+\cdots +{{a}_{k-1}}\xi
_{0}^{k-1}+{{a}_{k}}\xi _{0}^{k}}}.$$

Then we have

$${{a}_{0}}{{x_0}^{k}}+{{a}_{1}}{{x_0}^{k-1}}y_0+{{a}_{2}}{{x_0}^{k-2}}{{y_0}^{2}}+\cdots
+{{a}_{k-1}}x_0{{y_0}^{k-1}}+{{a}_{k}}{{y_0}^{k}}=$$
$$=a_{0}x^{k}_{0}+a_{1}x^{k}_{o}{\xi}_{0}+a_{2}x^{k}_{o}{\xi}^{2}_{0}+\cdots+a_{k-1}x^{k}_{o}{\xi}^{k-1}_{0}+a_{0}x^{k}{\xi}^{k}_{0}=$$

$$={\left(\frac{1}{\sqrt[k-1]{a_{0}+a_{1}{\xi}_{0}+a_{2}{\xi}^{2}_{0}+\cdots+a_{k-1}{\xi}^{k-1}_{0}+a_{k}{\xi}^{k}_{0}}}\right)^{k}}\cdot$$

$$\cdot({a_{0}+a_{1}{\xi}_{0}+a_{2}{\xi}^{2}_{0}+\cdots+a_{k-1}{\xi}^{k-1}_{0}+a_{k}{\xi}^{k}_{0}})=$$

$$=\frac{1}{\sqrt[k-1]{a_{0}+a_{1}{\xi}_{0}+a_{2}{\xi}^{2}_{0}+\cdots+a_{k-1}{\xi}^{k-1}_{0}+a_{k}{\xi}^{k}_{0}}}=x_{0}.$$

By suppose the number ${\xi}_{0}$ is root of the algebraical
equation (3.2), i.e.

$$ a_{k}{\xi}^{k+1}_{0}+\sum\limits_{i=0}^{k-1}{(a_{k-1-i}-b_{k-i})}{\xi}^{k-i}_{0}-b_{0}=0.$$

Hence we obtain the following equalities:

$$b_{0}+b_{1}{\xi}_{0}+b_{2}{\xi}^{2}_{0}+\cdots+b_{k}{\xi}^k_{0}=
a_{0}{\xi}_0+a_{1}{\xi}^2_{0}+a_{2}{\xi}^{3}_{0}+\cdots+a_{k}{\xi}^{k+1}_{0},$$

$$b_{0}+b_{1}{\xi}_{0}+b_{2}{\xi}^{2}_{0}+\cdots+b_{k}{\xi}^k_{0}=
{\xi}_0{a_{0}+a_{1}{\xi}_{0}+a_{2}{\xi}^{2}_{0}+\cdots+a_{k}{\xi}^{k}_{0}}.$$

Multiplying both sides of the last equality by the expression

$$\frac{1}{\left(a_{0}+a_{1}{\xi}_{0}+a_{2}{\xi}^{2}_{0}+\cdots+a_{k}{\xi}^{k}_{0}\right)^k}$$

we have

$${{b}_{0}}{{x}_{0}}^{k}+{{b}_{1}}{{x}_{0}}^{k-1}{{y}_{0}}+{{b}_{2}}{{x}_{0}}^{k-2}{{y}_{0}}^{2}+\cdots
+{{b}_{k-1}}{{x}_{0}}{{y}_{0}}^{k-1}+{{b}_{k}}{{y}_{0}}^{k}={{y}_{0}}.$$\\

This completes proof of the lemma 3.2.

\section{An application: Translational-invariant Gibbs measures for models on the Cayley tree $\Gamma^{k}$} \label{sec:intro}

A Cayley tree $\Gamma^k=(V,L)$ of order $k\geq 1$ is an infinite
homogeneous tree, i.e., a graph without cycles, with exactly $k+1$
edges incident to each vertices. Here $V$ is the set of vertices
and $L$ that of edges. Consider models where the spin takes values
in the set $[0,1]$, and is assigned to the vertices of the tree. A
configuration $\sigma$ on $V$ is defined as a function $x\in
V\mapsto\sigma (x)\in [0,1]$, the set of all configurations is
$[0,1]^V$. We consider the model $H$ on the $\Gamma^k$ by the
equality:
\begin{equation}\label{e1.1}
 H(\sigma)=-J\sum_{\langle x,y\rangle\in L}
\xi_{\sigma(x), \sigma(y)}, \,\ \sigma\in\Omega_{V}
\end{equation}
where $J \in R\setminus \{0\}$ and $\xi: (u,v)\in [0,1]^2\to
\xi_{uv}\in \mathbb{R}$ is a given bounded, measurable function.
As usually, $\langle x,y\rangle$ stands for the nearest neighbor
vertices.

Write $x<y$ if the path from $x^0$ to $y$ goes through $x$. Call
vertex $y$ a direct successor of $x$ if $y>x$ and $x,y$ are
nearest neighbors. Denote by $S(x)$ the set of direct successors
of $x$. Observe that any vertex $x\ne x^0$ has $k$ direct
successors and $x^0$ has $k+1$.

Let $h:\;x\in V\mapsto h_x=(h_{t,x}, t\in [0,1]) \in R^{[0,1]}$ be
mapping of $x\in V\setminus \{x^0\}$.

Now, we consider the following equation:
\begin{equation}\label{e5}
 f(t,x)=\prod_{y\in S(x)}{\int_0^1\exp(J\beta\xi_{tu})f(u,y)du \over \int_0^1\exp(J\beta{\xi_{0u}})f(u,y)du}.
 \end{equation}
Here, and below  $f(t,x)=\exp(h_{t,x}-h_{0,x}), \ t\in [0,1]$ and
$du=\lambda(du)$ is the Lebesgue measure.

 It is known that, for the
splitting Gibbs measure  for the model (\ref{e1.1}) to exist,  the
existence a solution of the equation (\ref{e5}) for any $x\in
V\setminus\{x^0\}$ is necessary and sufficient. Thus, we know that
the Gibbs measure $\mu$ for the model (\ref{e1.1}) depends  on the
function $f(t,x)$ and each the Gibbs measure corresponds to the
solution $f(t,x)$ of the equation (\ref{e5}).

A detailed definition of the splitting Gibbs measure for models
with nearest-neighbor interactions and the continuum set of spin
values on the Cayley tree $\Gamma^k$ is given in the works
\cite{er2010}, \cite{u2013}, \cite{g2011}. In the future, in place
term the splitting Gibbs measure, we will use the Gibbs measure.

Note, that the analysis of solutions to (\ref{e5}) is not easy.
It's difficult to give a full description for the given potential
function $\xi_{t,u}$. We study the Gibbs measures of the model
(\ref{e5}) in the case $f(t,x)=f(t)$ for all $x\in S(x)$. Such
Gibbs measure is called translation-invariant.

\begin{lemma}\label{l2.1.}\cite{ehr2012} Let $k\geq2$. The Hamiltonian
$H$ (\ref{e1.1}) has a translation-invariant Gibbs measure iff the
Hammerstein's operator $H_k$ has a positive eigenvalue, i.e. the
Hammerstein's equation
\begin{equation}\label{e2.2} H_{k}f=\lambda f, \,\ f\in C_{+}[0,1]
\end{equation}
has a nonzero positive solution for some $\lambda>0$.
\end{lemma}

 Moreover, if the number $\lambda_0>0$ is
an eigenvalue of the operator $H_{k}$, then an arbitrary positive
number is an eigenvalue of the operator $H_{k}$ (see Theorem 3.7
\cite{ehr2013}), where $k\geq2$. By the Lemma 4.1 and from the
proof of the Theorem 3.7 \cite{ehr2013} implies

\begin{lemma}\label{l2.2.} Let $k\geq 2$. A number $N^{tigm}(H)$ of
the translation-invariant Gibbs measures for the model
(\ref{e1.1}) the following equality holds:
$$N^{tigm}(H)=N_{+}^{fix}(H_k),$$ where $N_{+}^{fix}(T)$ is the
number of nontrivial positive fixed points of the operator $T$.
\end{lemma}

In this section we study translation-invariant Gibbs measures for
the wollowing model on the Cayley tree $\Gamma^k$:
$$
H(\sigma)=-\frac{1}{\beta} \sum \limits_{<x,y> \in L}\ln\left(a+ b
\sigma(x) \sigma(y)\right), \eqno{(4.4)}
$$
where there parameters $a$ and $b$ satisfy the conditions: $a>0,
\,\ b>0$ and $\beta=T^{-1}, T -$ temperature, $T>0.$

\textbf{Theorem 4.1.} The model (4.1) has unique
translational-invariant Gibbs measure for all $k\in \mathbb{N}$.

Proof. For the kernel $K(t,u)$ of the Hammerstain integral
operator $H_{k}$ we
 have

$${{\varphi }_{1}}\left( t \right)=1,\,\,\,\,{{\psi }_{1}}\left( t
\right)=a,\,\,\,\,{{\varphi }_{2}}\left( t\right)=t, \,\
\,\,\,\,{{\psi }_{2}}\left( t \right)=bt.$$

Thus, we need consider the following Hammerstein's operator

 $$(H_{k}f)(t)=\int\limits_0^1 (a+btu)f^{k}(u)du. \eqno{(4.5)}$$

Consequently, for $a_{i}$ and $b_{i}$ we get (see the section 2):

$$a_i= aC_k^{i}\int\limits_{0}^{1}u^{i}du=\frac{a}{i+1}\cdot\frac{k!}{(k-i)! \cdot i!}, \,\,\ i\in \{0,1,2, \cdot\cdot\cdot, k\};$$

$$b_i= bC_k^{i}\int\limits_{0}^{1}u^{i+1}du=\frac{b}{i+2}\cdot\frac{k!}{(k-i)! \cdot i!}, \,\,\ i\in \{0,1,2, \cdot\cdot\cdot, k\}.$$

Put

$$d_i=a_{i-1}-b_i=\frac{k!}{(k-i)! \cdot i!}\left( \frac{a}{k-i+1}  -\frac{b}{i+2} \right), \,\,\ , \,\,\ i\in \{1,2, \cdot\cdot\cdot, k\}.$$

 The inequality $d_i\geq 0, \,\,\ i\in \{1, 2, \cdot\cdot\cdot, k\}$ entails the following
 inequality:

 $$\frac{a}{k-i+1}  -\frac{b}{i+2}\geq0, \,\,\ i\in \{1,2, \cdot\cdot\cdot, k\}.$$

We define function

$$h(x)=\frac{a}{k-x+1}  -\frac{b}{x+2}$$

on $[1,k].$ Then we have

$$h'(x)=\frac{a}{(k-x+1)^{2}}  +\frac{b}{(x+2)^{2}}, \,\,\ x \in [1,k].$$

From here we get $h'(x)>0$ for all $x \in [1,k]$. Consequently,
the function $h(x)$ is increasing on the set $[1,k].$ Therefore,
for the numbers $d_1, d_2, \cdot\cdot\cdot, d_k$ the relations

$$d_1 \leq d_2 \leq \cdot\cdot\cdot \leq d_k$$

hold.  Then by the Theorem 2.2 the integral equation (2.1) has
unique positive solution. It means that, the model (4.4) has
unique translational-invariant Gibbs measure for all $k\in
\mathbb{N}$.

\end{document}